\documentclass[12pt,reqno]{amsart}
\usepackage[latin1]{inputenc}
\usepackage{amsmath}
\usepackage{amsfonts}
\usepackage{amssymb}
\usepackage{graphicx}
\usepackage{enumitem}
\usepackage{xcolor}
\usepackage{url}
\usepackage{amsthm}
\usepackage{bm}
\usepackage[left=3cm,right=3cm,top=2cm,bottom=2cm]{geometry}
{\theoremstyle{definition}
\newtheorem{definition}{Definition}[section]} 
\newtheorem{theorem}[definition]{Theorem}
\newenvironment{theorem*}[1]{{\bf Theorem #1} \begin{itshape}}{\end{itshape}}
\newtheorem{lemma}[definition]{Lemma}

\newenvironment{corollary*}[1]{{\bf Corollary #1} \begin{itshape}}{\end{itshape}}

\newenvironment{proposition*}[1]{{\bf Proposition #1} \begin{itshape}}{\end{itshape}}

\theoremstyle{remark}
\newtheorem{remark}[definition]{Remark}
\usepackage{geometry}
\newcommand{\N}{\mathbb{N}}

\newcommand{\Z}{\mathbb{Z}}
\newcommand{\Q}{\mathbb{Q}}
\newcommand{\R}{\mathbb{R}}

\newcommand{\Qs}{\mathbb{Q}_{S}}
\newcommand{\Zs}{\mathbb{Z}_{S}}

\newcommand{\Qv}{\mathbb{Q}_{\nu}}
\newcommand{\Zv}{\mathbb{Z}_{\nu}}

\newcommand{\cH}{\mathcal{H}}

\newcommand{\bx}{\mathbf{x}}
\newcommand{\by}{\mathbf{y}}
\newcommand{\ba}{\mathbf{a}}

\newcommand{\bb}{\mathbf{b}}
\newcommand{\bv}{\mathbf{v}}

\newcommand{\bt}{\boldsymbol{\tau}}

\newcommand{\lL}{(\ell)}

\newcommand{\supp}{\operatorname{supp}}

\title[Rectangular shrinking targets for actions on tori]{Rectangular shrinking targets for $\Z^m$ actions on tori: well and badly approximable systems}

\author{V. Beresnevich}
\address[V. Beresnevich]{Department of Mathematics, University of York, Heslington, York, YO10
5DD, United Kingdom}
\email{victor.beresnevich@york.ac.uk}

\author{S. Datta}
\address[S. Datta]{Department of Mathematics, University of Michigan, Ann Arbor, MI 48109, USA}
\email{shreyasi1992datta@gmail.com}

\author{A. Ghosh}
\address[A. Ghosh]{School of Mathematics, Tata Institute of Fundamental Research, Homi Bhabha Road, Colaba, Mumbai 400005, India}
\email{ghosh.anish@gmail.com}

\author{B. Ward}
\address[B. Ward]{La Trobe University Bendigo Campus, Edwards Rd, Flora Hill, Bendigo, Victoria 3552, Australia}
\email{Ben.Ward@latrobe.edu.au}

\thanks{A.\ G.\ gratefully acknowledges support from a grant from the Infosys foundation, a Department of Science and Technology, Government of India, Swarnajayanti fellowship and a grant from the Department of Atomic Energy, Government of India, under project $12-R\&D-TFR-5.01-0500$.\\
B.\ W.\ acknowledges support from Australian Research Council Discovery grant no. 200100994.}

\date{\today}
\begin{document}

\begin{abstract}
In this paper we investigate the shrinking target property for irrational rotations. This was first studied by Kurzweil (1951) and has received considerable interest of late. Using a new approach, we generalize results of Kim (2007) and Shapira (2013) by proving a weighted effective analogue of the shrinking target property. Furthermore, our results are established in the much wider $S$-arithmetic setting. 

\end{abstract}

\maketitle

\section{Introduction}\label{Intro}

The shrinking target problem for a dynamical system typically seeks to study visits of orbits of the dynamical system to some neighbourhoods of a point (regarded as `targets') whose sizes are `shrinking' over `time' as the orbits evolve. It is a manifestation of the chaotic behaviour of the dynamical system and is closely connected to questions in  Diophantine approximation. This paper is devoted to perhaps the oldest and simplest example of the shrinking target problem, namely the case of irrational rotations on the circle, or more generally, $\Z^m$ actions on tori. The `targets' in this case are the neighbourhoods of a point on the torus and the problem is to hit shrinking balls around this point by orbits of an irrational rotation. This problem is closely connected to \emph{inhomogeneous Diophantine approximation}. The subject can be said to have started with a result of Kurzweil \cite{Kurzweil} who provided a necessary and sufficient condition on the rotation angle to ensure that a Borel-Cantelli style `0-1' law holds in this setting. More recently, in \cite{MR2335077}, Kim proved the following beautiful result.
\begin{theorem}
Let $\theta$ be any irrational number. For almost every $s \in \R$ we have
$$\liminf_{n \to \infty} n \min_{b \in \Z} |n\theta - b - s| = 0.$$
\end{theorem}

Kim used techniques from the theory of continued fractions, specifically the Ostrowski expansion, and also obtained further refinements in \cite{MR3247068}. See also \cite{MR3346166, Kim2023} and \cite{Moshchevitin5} for recent developments in this direction. A higher dimensional version of Kim's result was proved by Shapira \cite{Shapira} using the tools of homogeneous dynamics, namely the classification of divergent orbits for diagonal flows on the space of lattices and considerations from the geometry of numbers. The irrationality condition, imposed by Kim in the one dimensional case, is replaced with a \emph{non-singularity} condition in higher dimensions, see Definition \ref{def:nonsingular} below. However, note that non-singularity coincides with irrationality in the one dimensional setting.  We will provide a bi-fold generalization of Kim's result as well as Shapira's result by firstly establishing the corresponding shrinking target property for rectangular targets (defined by weights) and secondly considering the much wider $S$-arithmetic setting. Our methods differ from those of Kim and Shapira; we rely on a softer approach using the technology of `constant invariant limsup sets' developed in our paper \cite{BDGW_null}.   

In recent years, the shrinking target problem for $\Z^m$ actions on tori has received considerable attention from both the dynamical systems as well as the number theory community. In addition to the works cited above, \cite{Bugeaud2003, MR1992153, FuchsKim, KimRamsWang, Kim2023} deal with Hausdorff dimension questions related to circle rotations and inhomogeneous Diophantine approximation. There have also been several developments on the  `badly approximable' problem beginning with \cite{BHKV}.

\subsection{Real tori: a special case}
We begin by stating a special case of our main theorem, namely a weighted generalization of Shapira's theorem.  
Let $\bm\tau=(\tau_1,\dots,\tau_n) \in \bm\R_{>0}^n$ and $\bm\eta =(\eta_1,\dots,\eta_m)  \in \bm\R_{>0}^m$ be  $n$- and $m$-tuples of non negative real numbers that will be referred to as {\em weights}. We will assume throughout that
$$
\sum_{i=1}^n\tau_i = \sum_{j=1}^m\eta_j\,.
$$
Given vectors $\bx=(x_1,\dots,x_n)\in\R^n$ and $\by=(y_1,\dots,y_m)\in\R^m$, we define the following quasi-norms associated with the weights:
$$
|\bx|_{\bm\tau}:=\max_{1\le i\le n}|x_i|^{1/\tau_i}\qquad\text{and}\qquad |\by|_{\bm\eta}:=\max_{1\le j\le m}|y_j|^{1/\eta_j}\,.
$$

\bigskip

In weighted Diophantine approximation, one replaces the usual supremum norms with weighted quasi-norms as above.  With this notation, we can now state a weighted analogue of Dirichlet's theorem. The proof follows from the usual argument using Minkowski's convex body theorem.\\

\noindent\textbf{Weighted Dirichlet's theorem:} {\em Let $\bm\tau$ and $\bm\eta$ be as above. Then for any real $m\times n$ matrix $\bm\alpha=(\alpha_{i,j})$ and any $H>1$ there exists $(\ba,\bb)\in\Z^{m+n}$ with $\ba\neq\bm0$ such that 
\begin{equation}
|\ba\bm\alpha-\bb|_{\bm\tau}<H^{-1}\qquad\text{and}\qquad |\ba|_{\bm\eta}\le H\,.
\end{equation}
}

\bigskip
Similarly, one can define a weighted analogue of singular vectors and non-singular vectors. 
\begin{definition}\label{def:nonsingular}
An $m\times n$ real matrix $\bm\alpha$ will be called {\em$(\bm\tau,\bm\eta)$-singular} if for any $\varepsilon>0$ there exists $H_0$ such that for all $H>H_0$ there exists $(\ba,\bb)\in\Z^{m+n}$ with $\ba\neq\bm0$ satisfying
\begin{equation}
|\ba\bm\alpha-\bb|_{\bm\tau}<\varepsilon H^{-1}\qquad\text{and}\qquad |\ba|_{\bm\eta}\le H\,.
\end{equation}
Otherwise, $\bm\alpha$ will be called {\em$(\bm\tau,\bm\eta)$-non-singular}. 
\end{definition}

We are now ready to state  the following weighted generalisation of the Kim-Shapira theorem.

\begin{theorem}\label{thm1}
Let $\bm\tau$ and $\bm\eta$ be as above and $\bm\alpha$ be any real $m\times n$ matrix which is $(\bm\tau,\bm\eta)$-non-singular. Then for almost every $\bm\gamma\in\R^n$ we have that for any $\varepsilon>0$ there are infinitely many $H\in\N$ such that
\begin{equation}
|\ba\bm\alpha-\bb-\bm\gamma|_{\bm\tau}<\varepsilon H^{-1}\qquad\text{and}\qquad |\ba|_{\bm\eta}\le H\,,
\end{equation}
or equivalently 
\begin{equation}
\liminf_{|\ba|_{\bm\tau}\to\infty}\; |\ba|_{\bm\eta}\min_{\bb\in\Z^n}|\ba\bm\alpha-\bb-\bm\gamma|_{\bm\tau}=0\,.
\end{equation}
\end{theorem}

In Section \ref{sec:sketch} we sketch a proof of Theorem \ref{thm1}. As mentioned earlier, our main theorem is considerably more general and includes the $S$-arithmetic setting. We now introduce this setup and state our main result.

\section{$S$-arithmetic setup and the main theorem}

Let $S$ be a collection of valuations of $\Q$ with cardinality $l\ge1$. Let $S^{*}=S \backslash \{ \infty \}$. Let $\Qs$ denote the set of $S$-arithmetic points, that is
\begin{equation*}
\Qs:=\prod_{\nu \in S} \Qv,
\end{equation*}
where $\Qv$ is the completion of $\Q$ with respect to the valuation $|\cdot|_{\nu}$. In particular, if $\nu\in S^*$ is a prime number $p$, then $\Q_\nu=\Q_p$ is the field of $p$-adic numbers, and if $\nu=\infty$, then $\Q_\nu=\R$. 
Let
\begin{equation*}
\Zs:=\prod_{\nu \in S} \Zv\,,
\end{equation*}
where for $\nu\in S^*$ we have that $\Z_\nu=\{x\in\Q_\nu:|x|_\nu\le1\}$ is the set of $\nu$-adic integers, and $\Z_\infty= [0,1)$. Let $\mu_{S}$ be the $S$-arithmetic Haar measure, normalised by $\mu_{S}(\Zs)=1$. Note that $\mu_{S}$ is simply the product of the measures $\mu_\nu$ over each $\Q_{\nu}$ normalised so that $\mu_\nu(\Z_\nu)=1$ for $\nu\in S^*$ with $\mu_\infty$ being Lebesgue measure over $\R=\Q_\infty$, {\em i.e.}
\begin{equation*}
\mu_{S}=\prod_{\nu \in S} \mu_{\nu}\,.
\end{equation*}
Similarly, in higher dimensional settings, associated with systems of linear  forms, let $\mu_{S,m\times n}$ denote the Haar measure on $\Qs^{mn}$, normalised by $\mu_{S,m\times n}(\Zs^{mn})=1$. We will denote any $\bx\in\Q_S^{mn}$ as $\bx=(\bx_\nu)=(\bx_\nu)_{\nu\in S}$, where $\bx_\nu=(x_{i,j,\nu})\in\Q_\nu^{mn}$, that is $(x_{i,j,\nu})$ is an $m\times n$ matrix over $\Q_\nu$ for each $\nu$. Let 
\begin{equation} \label{omega_definition}
\omega=\begin{cases}
m+n \quad \text{ if } \infty \not \in S, \\
\quad m \quad \, \, \, \,\, \text{ if } \infty \in S.
\end{cases}
\end{equation}

Let $\bt=(\tau_{i,\nu})_{1\le i\le n,\;\nu\in S}$ and $\bm\eta=(\eta_\ell)_{1\le\ell\le \omega}$ be two collections of positive real numbers, which will be referred to as {\em weights}, such that
 \begin{equation} \label{S_arithmetic_tau_condition}
 \sum_{i=1}^{n} \sum_{\nu\in S}\tau_{i,\nu}=\omega=\sum_{\ell=1}^\omega \eta_\ell\,.
 \end{equation}

 Given weights $\bt=(\tau_{i,\nu})_{1\leq i\leq n, \nu\in S}$, we define the $\bt$-semi-norm of $\bx\in\Q_S^n$ as
 $$\vert \bx\vert_{\bm\tau}=\max_{\nu\in S}\max_{1\leq i\leq n}\vert x_{i,\nu}\vert_\nu^{1/\tau_{i,\nu}}.$$
 Also, when $\infty\in S$, we will separate out the weights at the infinite and finite places of $S$ by introducing $\bm\tau_\infty:=(\tau_{i,\infty})$ and $\bm\tau^\star=(\tau_{i,\nu})_{\nu\in S^\star}.$ 
The following is  an $S$-arithmetic version of Dirichlet theorem with weights, see for example \cite[Corollary 6.3]{BDGW_null}.

\begin{theorem}\label{cor6.3}
Let $\bt=(\tau_{i,\nu})_{1\le i\le n,\;\nu\in S}$ and $\bm\eta=(\eta_\ell)_{1\le\ell\le \omega}$ be two collections weights satisfying \eqref{S_arithmetic_tau_condition}. Then for any $\bx=(x_{i,j,\nu})\in \Zs^{mn}$ and any $H>1$ there exists $(\ba,\bb)=(a_{1}, \dots , a_{m},b_1,\dots,b_n)\in\Z^{m+n}\setminus\{\bm0\}$ satisfying\begin{align}
\label{vb042}
\left|a_{1}x_{i,1,\nu}+ \dots + a_{m}x_{i,m,\nu}-b_{i} \right|_{\nu}&\le\nu H^{-\tau_{i,\nu}} 
&&(1 \leq i \leq n,\;
\nu \in S^*)\,\\
\label{vb043}
\left|a_{1}x_{i,1,\nu}+ \dots + a_{m}x_{i,m,\nu}-b_{i} \right|_{\nu}&<H^{-\tau_{i,\nu}} 
&&(1 \leq i \leq n,\;\nu=\infty)\quad\text{if $\infty\in S$\,}\\
\label{integer_condition3}
|a_{j}| &\leq H^{\eta_j} 
&&(1 \leq j \leq m)\,,\\
\label{integer_condition4}
|b_{i}| &\leq H^{\eta_{m+i}}  
&&(1 \leq i \leq n)\hspace{10.5ex}\text{if $\infty\not\in S$\,.}
\end{align}
\end{theorem}

\medskip

\begin{remark}
In the above theorem if $\infty\in S$ or 
\begin{equation}\label{vb044}
\eta_{m+i}<\sum_{\nu\in S}\tau_{i,\nu}\qquad\text{for all $1\le i\le n$}
\end{equation}
then for large enough $H$ for any $(\ba,\bb)\in\Z^{m+n}\setminus\{\bm0\}$, inequalities \eqref{vb042}---\eqref{integer_condition4} imply that $\ba\neq\bm0$.
\end{remark}

\medskip

\noindent Next, we define weighted singular and non-singular matrices in the $S$-arithmetic setting.

\begin{definition}
Let $\bm\alpha=(\alpha_{i,j,\nu})\in\Zs^{mn}$ be given. Then $\bm\alpha$ will be called {\em$(\bm\tau,\bm\eta)$-singular} if for any $\varepsilon>0$ there exists $H_0>0$ such that for all $H>H_0$ there exists $(\ba,\bb)\in\Z^{m+n}\setminus\{\bm0\}$ satisfying
\begin{equation}\label{eqn: singular}
|\ba\bm\alpha-\bb|_{\bm\tau}<\varepsilon H^{-1}\qquad\text{and}\qquad H\geq \left\{\begin{aligned}
    &|(\ba,\bb)|_{\bm\eta} \quad\text{if }\infty\notin S\\
    & |\ba|_{\bm\eta} \hspace{6.5ex}\text{if }\infty\in S.
    \end{aligned}\right.
\end{equation}
Otherwise, $\bm\alpha$ will be called {\em$(\bm\tau,\bm\eta)$-non-singular}. 
\end{definition} 

\medskip

\noindent Now we are fully ready to state our main Theorem:

\begin{theorem}\label{thm: S main theorem}
Let $\bm\tau$ and $\bm\eta$ be as above and $\bm\alpha$ be any $m\times n$ matrix with entries in $\Q_S$. Suppose $\bm\alpha$ is $(\bm\tau,\bm\eta)$-non-singular. Then for almost every $\bm\gamma\in\Q_S^n$ we have that for any $\varepsilon>0$ there are infinitely many $H\in\N$ such that for some $(\ba,\bb)\in\Z^{m+n}$ we have that
\begin{equation}\label{eqn: Dirichlet coro 1}
|\ba\bm\alpha-\bb-\bm\gamma|_{\bm\tau}<\varepsilon H^{-1}\qquad\text{and}\qquad H\geq \left\{\begin{aligned}
    &|(\ba,\bb)|_{\bm\eta} \quad\text{if }\infty\notin S\\
    & |\ba|_{\bm\eta} \hspace{6.5ex}\text{if }\infty\in S.
    \end{aligned}\right.\,.
\end{equation}
Equivalently, for almost every $\bm\gamma\in\Q_S^n$ we have that 
\begin{equation}\label{vb14}
\liminf_{\ba\in\Z^m}\; |\ba|_{\bm\eta}\min_{\bb\in\Z^n}|\ba\bm\alpha-\bb-\bm\gamma|_{\bm\tau}=0\qquad\text{if $\infty\in S,$}
\end{equation}
and
\begin{equation}\label{vb15}
\liminf_{(\ba,\bb)\in\Z^{m+n}}\; |(\ba,\bb)|_{\bm\eta}\;|\ba\bm\alpha-\bb-\bm\gamma|_{\bm\tau}=0\qquad\text{if $\infty\notin S.$}
\end{equation}
\end{theorem}
\medskip
 
\begin{remark}
If we have finitely many solutions $(\ba,\bb)\in\Z^{m+n}$ or $\ba\in\Z^m$ for $\infty\notin S$ and $\infty\in S$ respectively, to \eqref{eqn: Dirichlet coro 1} for infinitely many $H$, then $\ba\bm\alpha-\bb=\bm\gamma$ for some $(\ba,\bb)\in\Z^{m+n},$ and there are only countably many such $\bm\gamma.$ Thus, \eqref{vb14}/\eqref{vb15} follows immediately from \eqref{eqn: Dirichlet coro 1}. The converse in this equivalence is obtained even easier by setting $H$ to be $|(\ba,\bb)|$ or $|\ba|$ depending on the case we are in rounded to the next integer. 
\end{remark}

Prior to embarking onto the proof of the main theorem, in the next section, we sketch a proof of Theorem \ref{thm1}. It reveals the main ideas while avoiding some of the necessarily technical parts of the proof of Theorem \ref{thm: S main theorem}. In section \ref{sec:auxiliary} we introduce a key measure theoretic tool from our paper \cite{BDGW_null}. Section \ref{sec:proofmain} is devoted to the proof of the main Theorem.

\section{A sketch of a proof of Theorem~\ref{thm1}}\label{sec:sketch}

In what follows $\ll$ will denote the inequality $\le$ with an unspecified constant factor, and $\asymp$ will mean both $\ll$ and $\gg$. Let $\bm\alpha$ be $(\bm\tau,\bm\eta)$-non-singular. Then for arbitrarily large $t>0$ the first Minkowski minima of the lattice
$$
g_tu_{\bm\alpha}\Z^{m+n}
$$
is bounded below by a fixed constant, where
$$
g_t={\rm diag}\{e^{\tau_1t},\dots,e^{\tau_nt},e^{-\eta_1t},\dots,e^{-\eta_mt}\}
$$
and
$$
u_{\bm\alpha}=\left(\begin{array}{cc}
    \bm\alpha & I_m \\
     I_n & 0
\end{array}\right)\,.
$$
Fix one of these $t$. Then, by Minkowski's second theorem, the last minimum of $g_tu_{\bm\alpha}\Z^{m+n}$ is bounded above by a constant, say $C_0>0$ independent of $t$, and we have a basis of this lattice
$$
\bv_\ell=g_tu_{\bm\alpha}(\ba_\ell,\bb_\ell)^T\qquad(1\le \ell\le m+n)
$$
of vectors in the ball of radius $C_0$. Then, there are real numbers $x_1,\dots,x_{m+n}$ such that
$$
\sum_{\ell=1}^{m+n} x_\ell\bv_\ell=(\bm\gamma,\bm0)^T\,.
$$
Round each $x_\ell$ to the nearest non-zero integer $y_\ell$ and define
$$
(\ba,\bb)^T=\sum_{\ell=1}^{m+n} y_\ell (\ba_\ell,\bb_\ell)^T\in\Z^{m+n}\,.
$$
Then
$$
\left|g_tu_{\bm\alpha}(\ba,\bb)^T-(\bm\gamma,\bm0)^T\right|_\infty=\left|\sum_{\ell=1}^{m+n} y_\ell g_tu_{\bm\alpha}(\ba_\ell,\bb_\ell)^T-(\bm\gamma,\bm0)^T\right|_\infty=\left|\sum_{\ell=1}^{m+n} y_\ell \bv_\ell-(\bm\gamma,\bm0)^T\right|_\infty\ll 1
$$
which is the same as the system
\begin{align*}
    \left|\sum_{j=1}^ma_j\alpha_{i,j}-b_i-\gamma_i\right|\ll e^{-\tau_it}\qquad (1\le i\le n)\hspace*{1.6ex}\\    |a_j|\ll e^{\eta_jt}\qquad(1\le j\le m)\,.
\end{align*}
Setting $H\asymp e^t$ we get the following version of Dirichlet's theorem.

\bigskip

\noindent\textbf{Twisted weighted asymptotic Dirichlet for non-singular matrices:} {\em Let $\bm\tau$ and $\bm\eta$ be as above, and $\bm\alpha$ be any $(\bm\tau,\bm\eta)$-non-singular real $m\times n$ matrix. Then  there is a constant $C>0$ depending on $\bm\alpha$ only and an increasing unbounded sequence of $H_i>0$ such that for every $H=H_i$ and for every $\bm\gamma \in \R^{n}$ there exists $(\ba,\bb)\in\Z^{m+n}$ such that 
\begin{equation}\label{vb5}
|\ba\bm\alpha-\bb-\bm\gamma|_{\bm\tau}<CH^{-1}\qquad\text{and}\qquad |\ba|_{\bm\eta}\le H\,.
\end{equation}
Furthermore, for almost all $\bm\gamma$ there are infinitely many different $(\ba,\bb)$ arising from \eqref{vb5}.
}

\medskip

Theorem~\ref{thm1} now follows as a rather straightforward application of Lemma~\ref{CI_product} below (which in turn is Lemma 5.7 from \cite{BDGW_null}). The rather simple details can also be inferred from the proof of the same ilk in the $S$-arithmetic setting given in \S\ref{sec5.3}.

\section{An auxiliary result}\label{sec:auxiliary}

In this section we introduce an auxiliary result from \cite{BDGW_null}.
Fix an integer $n \geq 1$, and for each $1 \leq i \leq n$ let $(X_{i}, d_{i}, \mu_{i})$ be a metric space equipped with a $\sigma$-finite Borel regular measure $\mu_{i}$. Let $(X,d,\mu)$ be the product space with $X=\prod_{i=1}^{n}X_{i}$, $\mu=\mu_{1}\times\cdots\times\mu_n$ being the product measure, and
\begin{equation*}
d(\bx^{(1)},\bx^{(2)})=\max_{1 \leq i \leq n}d_{i}(x_i^{(1)},x_i^{(2)})\,, \qquad \text{where }\bx^{(j)}=(x_1^{(j)},\dots,x_n^{(j)})\,\,\text{for }j=1,2.
\end{equation*}

\begin{lemma}[\cite{BDGW_null} Lemma 5.7] \label{CI_product} 
Let $n\in\N$. For each $1\le j\le n$ let $(X_j,d_{j},\mu_j)$ be a metric measure space equipped with a $\sigma$-finite doubling Borel regular measure $\mu_{j}$. Let $X=\prod_{j=1}^{n}X_{j}$ be the corresponding product space, $d=\max_{1 \leq j \leq n}d_{j}$ be the corresponding metric, and $\mu=\prod_{j=1}^{n}\mu_{j}$ be the corresponding product measure. Let $(S_i)_{i\in\N}$ be a sequence of subsets of $\supp \mu$ and $(\bm\delta_{i})_{i \in \N}$ be a sequence of positive $n$-tuples $\bm\delta_{i}=(\delta^{(1)}_{i},\dots,\delta^{(n)}_{i})$ such that $\delta^{(j)}_{i} \to 0$ as $i \to \infty$ for each $1 \leq j \leq n$. Let
\begin{equation*}
\Delta_n(S_{i}, \bm\delta_{i})=\{\bx \in X: \, \exists \, \ba \in S_{i} \, \, \, d_{j}(a_{j}, x_{j}) < \delta_{i}^{(j)} \, \, \forall \, \, 1 \leq j \leq n\}\,,
\end{equation*}
where $\bx=(x_1,\dots,x_n)$ and $\ba=(a_1,\dots,a_n)$.
Then, for any $\bm C=(C_{1},\dots , C_{n})$ and $\bm c=(c_{1},\dots , c_{n})$ with $0<c_j\le C_{j}$ for each $1 \leq j \leq n$
\begin{equation}\label{vb25}
\mu\left( \limsup_{i \to \infty} \Delta_n(S_{i}, \bm C\bm \delta_{i}) \;\setminus\;\limsup_{i \to \infty} \Delta_n(S_{i}, \bm c\bm\delta_{i}) \right)=0\,,
\end{equation}
where $\bm c\bm\delta_j=(c_1\delta^{(1)}_i,\dots,c_n\delta^{(n)}_i)$ and similarly $\bm C\bm\delta_j=(C_1\delta^{(1)}_i,\dots,C_n\delta^{(n)}_i)$.
\end{lemma}

\bigskip

\section{ Proof of Theorem \ref{thm: S main theorem}}\label{sec:proofmain}
First, we prove the following theorem which is a twisted weighted asymptotic Dirichlet for non-singular matrices.

\begin{theorem}\label{thm: nonsingular Dirichlet coro}
Let $S$, $\bm\tau$ and $\bm\eta$ be as above, and $\bm\alpha$ be any $(\bm\tau,\bm\eta)$-non-singular  $m\times n$ matrix over $\Qs$. Then there is a constant $C>0$ depending on $\bm\alpha$ only and an increasing unbounded sequence $(H_i)_{i\in\N}$ of positive integers such that for every $i\in\N$ and any $\bm\gamma\in\Q_S^n$ there exists $(\ba,\bb)\in\Z^{m+n}\setminus\{\bm0\}$ such that 
\begin{equation}\label{eqn: inhomo Dirichlet for nonsingular}
|\ba\bm\alpha-\bb-\bm\gamma|_{\bm\tau}<CH_i^{-1}\qquad\text{and}\qquad H_i\geq \left\{\begin{aligned}
    &|(\ba,\bb)|_{\bm\eta} \quad\text{if }\infty\notin S\\
    & |\ba|_{\bm\eta} \hspace{6.5ex}\text{if }\infty\in S.
    \end{aligned}\right.
\end{equation}
\end{theorem}

\medskip

In what follows, given a lattice $\Gamma\subset\R^{m+n}$ and a convex body $B\subset\R^{m+n}$ symmetric about the origin, 
\begin{equation}\label{minimaoreder}
\lambda_1(\Gamma,B)\le \cdots \le \lambda_{m+n}(\Gamma,B)
\end{equation}
will denote the Minkowski minima of $\Gamma$ with respect to $B$, that is
$$
\lambda_i(\Gamma,B):=\inf\big\{\lambda>0:\mathrm{rank}(\Gamma\cap\lambda B)\ge i\big\}\,.
$$

\medskip

\subsection{ Proof of Theorem~\ref{thm: nonsingular Dirichlet coro} in the case $\infty\notin S$} 

Without loss of generality, we will assume that $\bm\alpha\in\Z_S^{mn}$ and $\bm\gamma\in \Z_S^n.$ For any given $\varepsilon>0$ and $H>1$
define 
\begin{equation}\label{Gamma}
\Gamma_{\bm\alpha}(\varepsilon, H):=\Big\{(\ba,\bb)\in\Z^{m+n}~:~|\ba\bm\alpha+\bb|_{\bm\tau}<\varepsilon H^{-1},\; \vert (\ba,\bb)\vert_\nu\leq \tfrac{1}{\nu} \quad\text{for all }\nu\in S\Big\}\,.
\end{equation}
It is readily seen, as a consequence of the strong triangle inequality for every $\nu\in S$, that $\Gamma_{\bm\alpha}(\varepsilon, H)$ is a sub-lattice of $\Z^{m+n}$. 
Furthermore, observe that 
\begin{equation}\label{VolumeOfLattice}
\text{Vol}\big(\R^{m+n}/\Gamma_{\bm\alpha}(\varepsilon, H)\big)\;\le\; \varepsilon^{-m-n}H^{\omega}\prod_{\nu\in S}\nu^{m+2n}\,.
\end{equation}
Next, let us consider the set 
$$K_H:=\Big\{(\ba,\bb)\in\R^{m+n}~:~\vert(\ba,\bb)\vert_{\bm\eta}\leq H\Big\}\,,
$$
which is obviously convex and symmetric about the origin. Observe that
\begin{equation}\label{VolumeOfBody}
\text{Vol}(K_H)=2^{m+n}H^\omega\,.
\end{equation}
Since $\bm\alpha$ is $(\bm\tau,\bm\eta)$-non-singular, there exists $\varepsilon>0$ and an unbounded subset $\cH$ of positive real numbers such that for any $H\in\cH$ system \eqref{eqn: singular} does not have any non-zero solution $(\ba,\bb)\in\Z^{m+n}$. 
This implies that
\begin{equation}\label{eqn: lower bound on first minima}
\lambda_1(\Gamma_{\bm\alpha}(\varepsilon, H), K_{H})\ge1\qquad\text{for all }H\in\cH\,.
\end{equation}
By Minkowski's second theorem on successive minima, we also have that
$$
\text{Vol}(K_{H})\prod_{\ell=1}^{m+n}\lambda_\ell(\Gamma_{\bm\alpha}(\varepsilon, H), K_{H})\leq 2^{m+n} \text{Vol}(\R^{m+n}/\Gamma_{\bm\alpha}(\varepsilon, H)).
$$
Hence, by Equations \eqref{minimaoreder}, \eqref{VolumeOfLattice}, \eqref{VolumeOfBody} and \eqref{eqn: lower bound on first minima}, for any $H\in\cH$ we have that 
\begin{equation}\label{upperbound1}
    \lambda_{m+n}(\Gamma_{\bm\alpha}(\varepsilon, H), K_{H}) \;\le\; C_0:=\varepsilon^{-m-n}\prod_{\nu\in S}\nu^{m+2n}\,.
\end{equation}
Fix any $H\in\cH$. Then, by \eqref{upperbound1}, there are $m+n$ linearly independent vectors $\{(\ba^{(\ell)},\bb^{(\ell)})\}_{\ell=1}^{m+n}$ in $\Gamma_{\bm\alpha}(\varepsilon, H)\cap C_0K_{H}$. Hence, by  \eqref{upperbound1} and the definitions of $K_H$ and $\Gamma_{\bm\alpha}(\varepsilon, H), K_{H})$, we have that
\begin{equation}\label{eqn: basis satisfying eqn}
\begin{aligned}
    &|\ba^{(\ell)}\bm\alpha+\bb^{(\ell)}|_{\bm\tau}<\varepsilon H^{-1}\\[0ex]
    & \vert a^{(\ell)}_j\vert \leq C_0H^{\eta_j}\hspace*{8ex} (1\le j\le m)\\[0ex]
    & \vert b^{(\ell)}_i\vert \leq C_0 H^{\eta_{m+i}}\hspace{6ex}(1\le i\le n)\\[0ex]
    &  \vert (\ba^{(\ell)},\bb^{(\ell)})\vert_\nu\leq \tfrac{1}{\nu}\hspace*{9.5ex}(\nu\in S).
    \end{aligned}
\end{equation}
Recall that for $H\in\cH$ the system \eqref{eqn: singular} does not have a non-zero integer solution. Therefore, for every $\ell=1,\dots,m+n$, we have that $\vert(\ba^{(\ell)},\bb^{(\ell)})\vert_{\bm\eta}>H$. This means that 
\begin{equation*}
\text{$\vert a_{j^\star}^{(\ell)}\vert >H^{\eta_{j^\star}}$ for some $1\le j^\star\le m$ \quad or \quad $\vert b_{i^\star}^{(\ell)}\vert >H^{\eta_{m+i^\star}}$ for some $1\le i^\star\le n$.}
\end{equation*}
Let us denote this $j^\star$ or $i^\star$, which depends on $\ell$, as $\ell^\star.$
So we can rewrite the above equation as follows, 
\begin{equation}\label{ijstar}
\vert b_{ \ell^\star}^{(\ell)}\vert >H^{\eta_{m+\ell^\star}}, 
\text{ or } 
\vert a_{\ell^\star}^{(\ell)}\vert >H^{\eta_{l^\star}}. 
\end{equation}

Now, fix the unique solution $(x^{(\ell)})_{\ell=1}^{m+n}\in\Q_S^{m+n}$ to the linear equation 
\begin{equation}\label{eqn: gamma as span}
(\mathbf{0},\bm\gamma)= \sum_{\ell=1}^{m+n} x^{(\ell)} (\ba^{(\ell)},\bb^{(\ell)})\,,
\end{equation}
which exists for the vectors $\{(\ba^{(\ell)},\bb^{(\ell)})\}_{\ell=1}^{m+n}$ are linearly independent.

Suppose that $\vert b^{(1)}_{1^\star}\vert >H^{\eta_{m+1^\star}}$, which is one of the two possibilities in \eqref{ijstar} for $\ell=1.$ See Remark \ref{remark:choosing b wlog} as to why we can assume this without loss of generality. Then for $\ell=1,\dots,m+n$, by the strong approximation theorem we choose $r^{(\ell)}\in\Q$ such that
\begin{equation}\label{eqn: strong approximation without infty}
\begin{aligned}
& \begin{cases}&|r^{(\ell)}-2\prod_{\nu\in S}\nu|\leq\prod_{\nu\in S}\nu \text{ if }  b^{(\ell)}_{1^\star}>0, \\
& |r^{(\ell)}+2\prod_{\nu\in S}\nu|<\prod_{\nu\in S}\nu, \text{ if } b^{(\ell)}_{1^\star}\leq 0,
\end{cases}\\[0ex]
&
|r^{(\ell)}-x^{(\ell)}|_S\leq 1 ,\\
& 
|r^{(\ell)}|_{\sigma}\leq 1
\quad \text{for all} \ \text{ primes }\sigma\notin S.
\end{aligned}
\end{equation}
Now, define 
\begin{equation}\label{ab}
\ba:= \sum_{\ell=1}^{m+n} r^{(\ell)}\ba^{(\ell)}\qquad\text{and}\qquad
 \bb:= \sum_{\ell=1}^{m+n} r^{(\ell)}\bb^{(\ell)}\,,
\end{equation}
which are thus rational vectors.
Then
$$
\ba\bm\alpha+\bb=\sum_{\ell=1}^{m+n}r^{(\ell)}(\ba^{(\ell)}\bm\alpha+\bb^{(\ell)}).
$$
First, we claim that $\ba\in\Z^{m}.$ 
From the last inequalities of \eqref{eqn: strong approximation without infty}, \eqref{ab} and the fact that $\ba^{(\ell)}$ and $\bb^{(\ell)}$ are integer vectors, we get that 
\begin{equation}\label{eqn: complement of S norm of q}\vert \ba\vert_\sigma\leq 1\quad\text{and}\quad\vert \bb\vert_\sigma\leq 1\qquad \text{for all} \ \text{ primes }\sigma\notin S.
\end{equation}
Note that, by \eqref{eqn: gamma as span},
$$
\sum_{\ell=1}^{m+n} x^{(\ell)} \ba^{(\ell)}=\mathbf{0}\,.
$$
Therefore, by the left hand side of \eqref{ab}, we get that
\begin{equation}\label{eqn: rewriting q}
\ba=\sum_{\ell=1}^{m+n}(r^{(\ell)}-x^{(\ell)})\ba^{(\ell)}.
\end{equation}
This, together with the second inequality in \eqref{eqn: strong approximation without infty}, implies that
$\vert \ba\vert_S\leq 1$.
Combining this with \eqref{eqn: complement of S norm of q} implies that $\ba\in\Z^{m}$.

\smallskip

Next, by \eqref{eqn: gamma as span} and \eqref{ab}, we have that 
$$\begin{aligned}
\ba\bm\alpha+\bb-\bm\gamma=&  \sum_{\ell=1}^{m+n}r^{(\ell)}(\ba^{(\ell)}\bm\alpha+\bb^{(\ell)})-  \sum_{\ell=1}^{m+n}x^{(\ell)}\bb^{(\ell)}-\sum_{\ell=1}^{m+n}x^{(\ell)}\ba^{(\ell)}\bm\alpha\\[1ex]
= & \sum_{\ell=1}^{m+n} (r^{(\ell)}-x^{(\ell)})(\ba^{(\ell)}\bm\alpha+\bb^{(\ell)}).
\end{aligned}
$$
Hence, using the first inequalities in \eqref{eqn: basis satisfying eqn} and the second inequalities in \eqref{eqn: strong approximation without infty}, we get from the above equation that
\begin{equation}\label{eqn: inhomo inequalities}
    \vert\ba\bm\alpha+\bb-\bm\gamma\vert_{\bm\tau}\leq \varepsilon H^{-1}.
\end{equation}
Now we write $\bb= \ba\bm\alpha+\bb-\bm\gamma+\bm\gamma -\ba\bm\alpha$. Using  Equation \eqref{eqn: inhomo inequalities} together with the fact that $\ba\in\Z^m$, $\bm\gamma\in\Z_S^n, \bm\alpha\in\Z_S^{mn},$ we conclude that $\bb\in\Z^n.$

Next, we claim that $(\ba,\bb)\neq \mathbf{0}.$
Note that $$b_{1^\star}= \sum_{\ell=1}^{m+n} r^{(\ell)} b^{(\ell)}_{1^\star}.$$
By the first inequalities in Equation \eqref{eqn: strong approximation without infty}, $r^{(\ell)}$ is positive if  $b^{(\ell)}_{1^\star}$ is positive and $r^{(\ell)}$ is non positive if $b^{(\ell)}_{1^\star}$ is so. Now, since $\vert b^{(1)}
_{1^\star}\vert >H^{\eta_{m+1^\star}},$ we are guaranteed that $\bb\neq 0$, which confirms our claim. 
Moreover, we get that $\vert b_{1^\star}\vert> H^{\eta_{m+1^\star}},$ which implies there exist infinitely many $(\ba,\bb)$ as solutions.

The last thing to show is that 
\begin{equation}\label{eqn: upper bound q,p}
\vert (\ba,\bb)\vert_{\bm\eta}\ll H.\end{equation}
By the second and third inequalities in Equation \eqref{eqn: basis satisfying eqn}, the first inequality in Equation \eqref{eqn: strong approximation without infty}, we get that
$$\vert a_j\vert \leq 3\prod_{\nu\in S}\nu (n+m)\varepsilon^{-m-n} \prod_{\nu\in S}\nu^{m+n} H^{\eta_k}, 1\leq j\leq m.$$
Similarly, we get 
$$\vert b_i\vert \leq 3\prod_{\nu\in S}\nu (n+m)\varepsilon^{-m-n} \prod_{\nu\in S}\nu^{m+n} H^{\eta_{m+i}}, 1\leq i\leq n.$$

Thus in the view of Equation \eqref{eqn: inhomo inequalities}, Equation \eqref{eqn: upper bound q,p}, and the discussion from above, the proof is complete.

\begin{remark}\label{remark:choosing b wlog}
    We can assume $\vert b^{(1)}_{1^\star}\vert>H^{\eta_{m+1^\star}}$ without loss of generality. If instead we had $\vert a^{(1)}_{1^\star}\vert > H^{\eta_{1^\star}}$, we will choose $r^{\lL}$ in Equation \eqref{eqn: strong approximation without infty} accordingly. Then a similar process as in the proof above will imply $\ba\neq \mathbf{0}.$
\end{remark}

\medskip

\subsection{Proof of Theorem~\ref{thm: nonsingular Dirichlet coro} in the case $\infty \in S$}. Let us denote $S^\star=S\setminus \infty$. Without loss of generality, we assume $\bm\alpha=(\bm\alpha_\star,\bm\alpha_\infty)\in{\Z^{mn}_{S^{\star}}}\times [0,1)^{mn}$ and $\bm\gamma=(\bm\gamma_\star,\bm\gamma_\infty)\in \Z_{S^{\star}}^{n}\times [0,1)^{n}.$
Let us define the lattice
$$\Gamma_{\bm\alpha}(\varepsilon, H):=\{(\ba,\bb)\in\Z^{m+n}~:~|\ba\bm\alpha_\star+\bb|_{\bm\tau^{\star}}<\varepsilon H^{-1}, \vert \ba\vert_\nu\leq \frac{1}{\nu} ~~\forall\nu\in S^\star\}.
$$
Next, let us consider the convex set 
$$K_H:=\{(\ba,\bb)\in\R^{m+n}~:~\vert \ba\bm\alpha_\infty+\bb\vert_{\bm\tau_\infty}<\varepsilon H^{-1}, \vert\ba\vert_{\bm\eta}\leq H\}.$$
If $\bm\alpha$ is {\em$(\bm\tau,\bm\eta)$-non-singular} then there exists $\varepsilon>0$ and an unbounded subset $\cH$ of positive real numbers such that for any $H\in\cH$ system
system \eqref{eqn: singular} does not have any non-zero solution. 
This implies that, 

\begin{equation}\label{eqn: lower bound on first minima infty}
\lambda_1(\Gamma_{\bm\alpha}(\varepsilon, H), K_{H})>1 \text{ for all } H\in \cH.
\end{equation}

Note that $$\text{Vol}(\R^{m+n}/\Gamma_{\bm\alpha}(\varepsilon, H))\leq \varepsilon^{-(l-1)n}H^{m-\sum_{k=1}^n\tau_{k,\infty}}\prod_{\nu\in S}\nu^{m},$$
and
$$
\text{Vol}(K_{H})=2^n\varepsilon^{n}H^{-\sum_{k=1}^n\tau_{k,\infty}} 2^m H^{m}.$$

Hence by Equation \eqref{eqn: lower bound on first minima infty}, we have for any $H\in \cH,$
\begin{equation}
    \lambda_{m+n}(\Gamma_{\bm\alpha}(\varepsilon, H), K_{H})\leq \varepsilon^{-ln} \prod_{\nu\in S}\nu^{m}.
\end{equation} Fix any $H\in \cH.$
We get $m+n$ many linearly independent $\{(\ba^{(\ell)},\bb^{(\ell)})\}_{l=1}^{m+n}\in\Gamma_{\bm\alpha}(\varepsilon, H)$ such that 
\begin{equation}\label{eqn: basis satisfying eqn infty}
\begin{aligned}
    &|\ba^{(\ell)}\bm\alpha_\star+\bb^{(\ell)}|_{\bm\tau^\star}<\varepsilon H^{-1},\\
    &|\ba^{(\ell)}\bm\alpha_\infty+\bb^{(\ell)}|_{\bm\tau_\infty}< \varepsilon^{-ln} \prod_{\nu\in S^\star}\nu^{m}\varepsilon H^{-1},\\
    & \vert a^{(\ell)}_k\vert \leq \varepsilon^{-ln} \prod_{\nu\in S^\star}\nu^{m} H^{\eta_k}, 1\leq k\leq m\\
    &  \vert \ba^{(\ell)}\vert_\nu\leq \frac{1}{\nu} ~~\forall\nu\in S^\star.
    \end{aligned}
\end{equation}
Since $\{(\ba^{(\ell)},\bb^{(\ell)})\}$ forms a basis in $\Q_S^{m+n},$
there exists a unique solution $\bx=(x^{\ell})_{l=1}^{m+n}\in\Q_S^{m+n}$ to the following system. 
\begin{equation}\label{eqn: gamma as span infty}
(D,\underbrace{\mathbf{0}}_{m-1},\bm\gamma)= (x^{(1)},\cdots,x^{(m+n)})\begin{bmatrix}
   & \ba^{(
   1)} &\cdots &\ba^{(m+n)}\\
   & \ba^{(1)}\alpha+\bb^{(1)} & \cdots &\ba^{(m+n)}\alpha+\bb^{(m+n)}
\end{bmatrix}^T,
\end{equation}
where $D\in\Q_S, D_\nu=0, \forall\nu\in S^\star $, and $D_\infty=H^{\eta_1}+\prod_{\nu\in S^\star}\nu \sum_{l=1}^{m+n} \vert a^{(\ell)}_{1}\vert.$
Hence, we have $$\bm\gamma= \sum_{\ell=1}^{m+n} x^{(\ell)} (\ba^{(\ell)}\alpha+\bb^{\ell}),~~~ \mathbf{0}= \sum_{\ell=1}^{m+m} x^{(\ell)}_{\nu}\ba^{(\ell)}, \forall~ \nu\in S^\star,$$
and 
$D_\infty = \sum_{\ell=1}^{m+n}x^{(\ell)}_{\infty} a^{\lL}_{1}.$

Using the strong approximation theorem, for $\ell=1,\cdots,m+n$ we get $r^{\lL}\in\Q$ such that
\begin{equation}\label{eqn: strong approximation with infty}
\begin{aligned}
& \vert r^{\lL}_{\infty}-x^{\lL}_{\infty}\vert \leq \prod_{\nu\in S^\star} \nu,\\
&|r^{\lL}-x^{\lL}|_{S^{\star}}\leq 1 ,\\
& |r^{\lL}|_{\nu}\leq 1
\quad \text{for all} \ \text{ primes }\nu\notin S.
\end{aligned}
\end{equation}

Now, let us take 
$ \ba:= \sum_{\ell=1}^{m+n} r^{\lL}\ba^{\lL},~~
 \bb:= \sum_{\ell=1}^{m+n}r^{\lL}\bb^{\lL}.$

Then
$$
\ba\bm\alpha+\bb=\sum_{\ell=1}^{m+n}r^{\lL}(\ba^{\lL}\bm\alpha+\bb^{\lL}).
$$
First, we claim that $\ba\in\Z^{m}.$ 
From the last inequalities in Equation \eqref{eqn: strong approximation with infty}, we get 
\begin{equation}\label{eqn: complement of S norm of q infty}\vert \ba\vert_\nu\leq 1, \vert \bb\vert_\nu\leq 1 ~\forall \nu\notin S.\end{equation}
Also by Equation \eqref{eqn: gamma as span infty},
 $\ba=\sum_{\ell=1}^{m+n} (r^{\lL}-x^{\lL}_{\nu})\ba^{\lL}$, $\nu\in S^\star$. 
This, together with the second inequality in Equation \eqref{eqn: strong approximation with infty}, implies
\begin{equation}\label{eqn: S norm of q infty}
\vert \ba\vert_{S^\star}\leq 1.
\end{equation} 
Hence by Equation \eqref{eqn: complement of S norm of q infty}, and Equation \eqref{eqn: S norm of q infty}, we conclude $\ba\in\Z^{m}.$

Since
$
\ba\bm\alpha+\bb-\bm\gamma=\sum_{\ell=1}^{m+n} (r^{\lL}-x^{\lL})(\ba^{\lL}\bm\alpha+\bb^{\lL}),
$
using the first inequality in Equation \eqref{eqn: basis satisfying eqn infty} and the second inequality in Equation \eqref{eqn: strong approximation with infty}, we get 

\begin{equation}\label{eqn: inhomo inequalities infty}
    \vert\ba\bm\alpha_\star+\bb-\bm\gamma_\star\vert_{\bm\tau^\star}\leq \varepsilon H^{-1}.
\end{equation}
Also, the second inequality in Equation \eqref{eqn: basis satisfying eqn infty} and the first inequality in Equation \eqref{eqn: strong approximation with infty}, we get
\begin{equation}\label{eqn: inhomo inequalities infty infty}
    \vert\ba\bm\alpha_\infty+\bb-\bm\gamma_\infty\vert_{\bm\tau_\infty}\leq \varepsilon^{-ln} 
    \prod_{\nu\in S^\star}\nu^{m}\varepsilon 
    H^{-1} {\left((n+m)\prod_{\nu\in S^\star}\nu\right)}^{\frac{1}{\tau'_\infty}},
\end{equation}
where $\tau'_\infty=\min_{i=1}^n \tau_{i,\infty}.$
Since $\bb= \ba\bm\alpha+\bb-\bm\gamma+\bm\gamma -\ba\bm\alpha$, using  Equation \eqref{eqn: inhomo inequalities infty} together with the fact that $\ba\in\Z^m$, $\bm\gamma_\star\in\Z_{S^\star}^n, \bm\alpha_\star\in\Z_{S^\star}^{mn},$ we conclude that $\bb\in\Z^n.$

Next, we claim that 
\begin{equation}\label{eqn: upper bound q infty}
\vert \ba\vert_{\bm\eta}\ll H.
\end{equation}

Since $a_k= \sum_{\ell=1}^{m+n}(r^{\lL}-x^{\lL}_{\infty})a^{\lL}_{k}$ for $k=2,\cdots,n,$ by the third inequality in Equation \eqref{eqn: basis satisfying eqn infty} and the first inequality in Equation \eqref{eqn: strong approximation with infty}, we get that
$$\vert a_k\vert \leq (n+m)\varepsilon^{-ln} \prod_{\nu\in S^\star}\nu^{m+1} H^{\eta_k}, 2\leq k\leq m.$$
Also, $$\begin{aligned} &a_1= \sum_{\ell=1}^{m+n} r^{\lL} a^{\lL}_{1}\implies a_1-D_\infty= \sum_{\ell=1}^{m+n} (r^{\lL}-x^{\lL}_{\infty}) a^{\lL}_{1}\\
& \implies \vert a_1\vert\leq  (n+m) \varepsilon^{-ln} \prod_{\nu\in S^\star}\nu^{m+1} H^{\eta_1}+D_\infty\ll H^{\eta_1}. \end{aligned}$$ 
Note that 
$$\begin{aligned}
& D_\infty= \sum_{\ell=1}^{m+n} (x^{\lL}_{\infty}-r^{\lL})a^{\lL}_{1} +\sum_{\ell=1}^{m+n} r^{\lL} a^{\lL}_{1}\\
&\implies H^{\eta_1}+ \prod_{\nu\in S^\star}\nu \sum_{\ell=1}^{m+n}\vert a^{\lL}_{1}\vert\leq \prod_{\nu\in S^\star}\nu \sum_{\ell=1}^{m+n}\vert a^{\lL}_{1}\vert+ \vert \sum_{\ell=1}^{m+n} r^{\lL} a^{\lL}_{1}\vert\\
& \implies H^{\eta_1} \leq \vert \sum_{\ell=1}^{m+n} r^{\lL} a^{\lL}_{1}\vert=\vert a_1\vert. \end{aligned}
$$
Thus, there are infinitely many $\ba\neq \mathbf{0}$ as solutions.
Therefore, in the view of Equation \eqref{eqn: inhomo inequalities infty}, Equation \eqref{eqn: inhomo inequalities infty infty}, and Equation \eqref{eqn: upper bound q infty} the proof is complete.

\subsection{ Proof of Theorem \ref{thm: S main theorem} using Theorem \ref{thm: nonsingular Dirichlet coro}}\label{sec5.3}

Let \begin{equation} \label{def: phi}
\phi(\ba,\bb):=\begin{cases}
\ba \quad \text{ if } \infty \not \in S, \\
 (\ba,\bb) \quad \, \, \, \,\, \text{ if } \infty \in S.
\end{cases}
\end{equation}
For any $\delta>0,$ let us define $$ B_{\phi(\ba,\bb)}(\delta, \bm\tau):=\left\{\bm\gamma\in\Z_S^n~\left|~ \begin{cases}\vert \ba\bm\alpha+\bb+\bm\gamma\vert_{\bm\tau} \quad \text{ if } \infty\not \in S,\\
\min_{\bb\in\Z^n}\vert \ba\bm\alpha+\bb+\bm\gamma\vert_{\bm\tau} \quad \text{ if } \infty\in S\end{cases}\leq \frac{\delta}{\vert \phi(\ba,\bb)\vert_{\bm\eta}}\right.
\right\}.$$
Then by Theorem \ref{thm: nonsingular Dirichlet coro}, $$\mu_{S,n}\left(\bigcup_{\delta>0} \limsup_{\vert \phi(\ba,\bb)\vert_{\bm\eta}\to\infty}B_{\phi(\ba,\bb)}(\delta, \bm\tau)\right)= \mu_{S,n} (\Z_S^n).$$
Hence by continuity of measure,
$$\lim_{\delta\to\infty+}\mu_{S,n}\left( \limsup_{\vert \phi(\ba,\bb)\vert_{\bm\eta}\to\infty}B_{\phi(\ba,\bb)}(\delta, \bm\tau)\right)= \mu_{S,n} (\Z_S^n).$$
By Lemma \ref{CI_product}, for any $\delta_1,\delta_2>0$ we have 
$$\mu_{S,n}\left( \limsup_{\vert \phi(\ba,\bb)\vert_{\bm\eta}\to\infty}B_{\phi(\ba,\bb)}(\delta_1, \bm\tau)\right)=\mu_{S,n}\left( \limsup_{\vert \phi(\ba,\bb)\vert_{\bm\eta}\to\infty}B_{\phi(\ba,\bb)}(\delta_2, \bm\tau)\right).$$
Hence, we get that
\begin{equation}
\mu_{S,n}\left(\bigcap_{\delta>0}\limsup_{\vert \phi(\ba,\bb)\vert_{\bm\eta}\to\infty}B_{\phi(\ba,\bb)}(\delta, \bm\tau)\right)= \mu_{S,n} (\Z_S^n)
\end{equation}
as required.



\bibliographystyle{plain}
\bibliography{badnull}
\end{document}